\documentclass{amsart}
\usepackage{a4}
\usepackage[french]{babel}
\usepackage[latin1]{inputenc}
\usepackage{xypic}
\xyoption{all}
\usepackage{epsfig}
\newtheorem{thm}{Th\'eor\`eme}[section]
 
\newtheorem{lem}[thm]{Lemme}
  
\newtheorem{cor}[thm]{Corollaire}
\newtheorem{prop}[thm]{Proposition}

\theoremstyle{definition}
\newtheorem{defn}[thm]{D\'efinition}
\newtheorem{exmp}[thm]{Exemple}

\theoremstyle{remark}
\newtheorem{rem}[thm]{Remarque} 

\newcommand {\spec}{\text{Spec }}
\newcommand {\spf}{\text{Spf }}
\newcommand {\ord}{\text{ord}}
\newcommand {\aut}{\text{Aut}}
\newcommand {\irr}{\text{Irr}}

\newcommand {\ind}{\text{Ind}}

\newcommand {\disc}{R[[Z]]}
\newcommand {\kuro}{\displaystyle{\frac{R[[Z_1,Z_2]]}{(Z_1Z_2-\pi^e)}}}
\newcommand {\kouro}{\mathcal A_e}
\begin{document}
\title[Rel\`evement galoisien...]{Rel\`evement galoisien des rev\^etements\\de courbes nodales}
\author{Yannick Henrio\footnotemark[1]}
\subjclass{Primary 14G20, 14L27 ; Secondary 14D15, 14E22}
\maketitle
\footnotetext[1]{Math\'ematiques Pures de Bordeaux, UPRES-A 5467 CNRS, 
Universit\'e Bordeaux I\\
\indent \indent 351 cours de la lib\'eration 33 405 Talence cedex, France\\
\indent \indent {\texttt e-mail : henrio@math.u-bordeaux.fr}}
\begin{abstract}
Let $R$ be a complete discrete valuation ring of mixed characteristics, with 
algebraically closed residue field $k$. We study the existence problem 
of equivariant liftings to $R$ of Galois covers of nodal curves over $k$. 
Using formal geometry, we show that this problem is actually a local one. 
We apply this local-to-global principle to obtain new results concerning 
the existence of such liftings.
\end{abstract}
\par
On fixe un corps alg\'ebriquement clos $k$ de 
caract\'eristique $p>0$, et $R$ un anneau de valuation discr\`ete complet 
de corps r\'esiduel $k$. On note $K$ le corps des fractions de $R$, 
de caract\'eristique nulle. Soit $Y_0$ une courbe alg\'ebrique 
projective sur $k$ nodale (c'est-\`a-dire connexe, r\'eduite, avec pour uniques 
singularit\'es des points doubles ordinaires), nous 
appellerons mod\`ele de $Y_0$ sur $R$ un couple 
$(Y,\psi)$, o\`u $Y$ est un sch\'ema normal, 
propre et plat sur $R$, de fibre g\'en\'erique lisse et 
g\'eom\'etriquement connexe sur $K$, et 
$\psi: Y \times_R k \to Y_0$ est un isomorphisme de $k$-sch\'emas.
Dans ce travail, nous \'etudions la question suivante :
\par
Si $G$ est un groupe fini de 
$k$-automorphismes de $Y_0$, agissant librement sur un ouvert dense, 
existe-t'il un triplet $(Y,\psi,\rho)$, o\`u :
\begin{itemize}
\item
$(Y,\psi)$ est un mod\`ele de $Y_0$ sur $R$.
\item
$\rho:G \to \aut_R Y$ est un homomorphisme injectif 
qui fait commuter le diagramme :
$$
  \diagram
     G \rto^{\rho}  \drto & \aut_R Y \dto\\
                   & \aut_k Y_0
  \enddiagram 
$$
(o\`u la fl\`eche verticale associe \`a un automorphisme $\sigma$ de 
$\aut_R Y$ l'automorphisme 
$\psi \circ \sigma_{|Y \times_R k}\circ \psi^{-1}$.)
\end{itemize}

Soit $(Y,\psi,\rho)$ un tel rel\`evement, le quotient 
$X:=Y/G$ est une courbe propre sur $R$, de fibre g\'en\'erique lisse sur $K$, 
et de fibre sp\'eciale nodale $X_0:=Y_0/G$. Ainsi le morphisme quotient 
$f:Y \to X$ est un rel\`evement $G$-galoisien du 
rev\^etement $f_0:Y_0 \to X_0$.
\par
Si le morphisme quotient $f_0:Y_0 \to X_0$ est \'etale, et si $(X, \phi)$ 
est un mod\`ele de $X_0$ sur $R$, il r\'esulte de la th\'eorie du groupe fondamental 
de Grothendieck qu'il existe un unique rev\^etement \'etale $G$-galoisien 
$Y \to X$ de fibre sp\'eciale $f_0$. Les obstructions 
au rel\`evement sont donc li\'es \`a  la ramification du morphisme $f_0$. 
En fait, si $Y_0$ (et donc $X_0$) est une courbe lisse, le r\'esultat 
pr\'ec\'edant s'\'etend au cas o\`u $f_0$ est mod\'er\'ement ramifi\'e. 
Toutefois, la ramification mod\'er\'ee impose des obstructions au 
rel\`evement local des points doubles : On doit supposer l'action de $G$ 
kumm\'erienne (voir \cite{Sa1}, th\'eor\`eme 5.7).
\par
Si $f_0$ est sauvagement ramifi\'e, il est 
n\'ecessaire de faire des hypoth\`eses sur les sous-groupes d'inertie pour 
obtenir des \'enonc\'es de rel\`evement. Par exemple, B. Green et M. Matignon 
ont montr\'e que si la courbe $Y_0$ est lisse, si $(X, \phi)$ est un 
mod\`ele de $X_0$ sur $R$, et si les groupes d'inertie 
sont tous cycliques de $p$-exposant inf\'erieur ou \'egal \`a $2$, il existe 
un rev\^etement $f:Y \to X$ galoisien de groupe $G$ qui rel\`eve $f_0$. 
En revanche, il n'y a plus unicit\'e du 
rel\`evement lorsque $X$ est fix\'e. Par ailleurs, une conjecture 
due \`a Oort (\cite{O1}, \cite{O2}) dit que si la courbe $Y_0$ est lisse et 
si les groupes d'inertie sont tous cycliques, il existe un rel\`evement sur un $R$ convenable. Dans ce travail, nous montrons le r\'esultat suivant:
\newline\newline
\noindent
{\bf Th\'eor\`eme}
{\it
Supposons l'action de $G$ kumm\'erienne et que pour un point 
ferm\'e $y$ de $Y_0$ :\\
1. Si y est un point lisse, le groupe d'inertie $I_y$ de $y$ est 
cyclique d'ordre 
$n(y)p^{r(y)}$, avec $(n(y),p)=1$ et $0\le r(y)\le 2$.\\
2. Si y est un point double, on a alors (i) ou bien (ii) :
\begin{itemize}
\item
(i) Le groupe d'inertie $I_y$ de $y$ est cyclique d'ordre $n(y)p^{r(y)}$ 
et l'action de $I_y$ sur les branches est triviale.
\item
(ii) Le groupe d'inertie $I_y$ de $y$ est di\'edral d'ordre $2n(y)p^{r(y)}$, 
avec $(n(y),p)=1$ et $0\le r(y)\le 2$, de pr\'esentation 
$$I_y=\langle \sigma,\tau|\sigma^{n(y)p^{r(y)}}=1,\tau^2=1,
\tau \sigma \tau= \sigma^{-1}\rangle.$$
Les branches du point double $y$ sont permut\'ees par $\tau$ et $\sigma$ 
induit un automorphisme d'ordre $n(y)p^{r(y)}$ de chacune des deux branches.
\end{itemize}
Alors, quitte \`a faire une extension finie de $K$, il existe un rel\`evement 
de $(Y_0,G)$.
}
\newline
\par
Green et Matignon ont \'etudi\'e le cas des courbes lisses \`a l'aide des 
m\'ethodes de la g\'eom\'etrie rigide. Nous utilisons ici la g\'eom\'etrie 
formelle, qui est sans doute mieux adapt\'ee aux questions 
concernant les mod\`eles entiers (on pourra consulter \cite{Ra} pour 
la comparaison des deux th\'eories). 
Plus pr\'ecis\'ement, nous nous sommes amplement inspir\'e du 
recollement formel (formal patching) \`a la Harbater (voir 
\cite{H} et \cite{H-S}). On en d\'eduit un principe local-global formel, 
qui montre que le probl\`eme de rel\`evement galoisien est essentiellement 
de nature locale. On est ainsi ramen\'e \`a construire des actions de groupes 
sur les disques et les couronnes formels.
\par
La premi\`ere partie consiste en quelques rappels sur la g\'eom\'etrie des 
disques et des couronnes formels. La seconde pr\'esente les m\'ethodes de 
recollement formel. En application, on d\'eduit l'existence d'un mod\`ele 
\`a \'epaisseurs fix\'ees sur $R$ pour 
toute courbe nodale projective. Nous montrons ensuite le principe 
local-global formel. Dans une troisi\`eme partie, nous construisons des 
rel\`evements locaux pour certaines actions de groupe sur un point double. 
La quatri\`eme partie est consacr\'ee \`a la d\'emonstration du th\'eor\`eme 
de rel\`evement \'enonc\'e plus haut.
\par
Ce travail constitue une partie des r\'esultats de ma th\`ese (\cite{Heth}), 
dont certains ont \'et\'e annonc\'es dans (\cite{He}). Dans un 
article ult\'erieur, nous \'etudierons la structure des automorphismes 
d'ordre $p$ du disque formel sur $R$ (\`a la suite de \cite{G-M 2}) et 
nous montrerons comment la g\'eom\'etrie formelle 
permet \'egalement de construire des automorphismes de disques ou de 
couronnes formels.
\newpage
\par
Nous utiliserons constamment les notations suivantes :
\begin{itemize}
\item
$\pi$ d\'esigne une uniformisante de $R$.
\item
Si $M$ est un $R$-module, on note $\overline{M}$ le 
$k$-espace vectoriel $M\otimes_R k$.
\item
Si $A$ est une $R$-alg\`ebre, $A_K$ (resp. $\bar A$) d\'esigne 
$A \otimes_R K$ (resp. $A \otimes_R k$), et pour un id\'eal premier $\frak p$ 
de $A$, on note $\hat A_{\frak p}=A_{\frak p}^{\wedge}$ le compl\'et\'e de 
$A_{\frak p}$ pour la topologie $\frak p$-adique.
\item
Si $T_1,\dots ,T_n$ sont des ind\'etermin\'ees, on note 
$R\{T_1,\dots ,T_n\}$ la $R$-alg\`ebre des s\'eries restreintes, 
c'est-\`a-dire la sous-$R$-alg\`ebre de $R[[T_1,\dots ,T_n]]$ 
form\'ee des s\'eries 
$$\sum_{(\nu_1,\dots,\nu_n)\in \mathbb{N}^n}
a_{\nu_1,\dots,\nu_n}T_1^{\nu_1}\dots T_n^{\nu_n}
\text{ telles que } 
\lim_{\nu_1+\dots+\nu_n \to +\infty} a_{\nu_1,\dots,\nu_n} =0.$$
\item
On d\'esigne par $R[[T]]\{T^{-1}\}$ la $R$-alg\`ebre des s\'eries de 
Laurent
\text{$f:=\sum_{\nu \in \mathbb{Z}}a_{\nu}T^{\nu}$,} \`a coefficients 
dans $R$, avec $\lim_{\nu \to -\infty}a_{\nu}=0$. 
C'est un anneau de valuation discr\`ete complet, d'uniformisante 
$\pi$ et de corps r\'esiduel $k((t))$.
\end{itemize}

\section{Disques et couronnes formels}
\subsection{Quelques d\'efinitions et rappels}
\begin{defn}
Soit $X$ un $R$-sch\'ema et $x$ un point ferm\'e de la fibre sp\'eciale de 
$X$, on appelle 
fibre formelle de $X$ au point $x$ le $R$-sch\'ema affine 
$\mathcal F(X,x):=\spec \hat \mathcal O_{X,x}$, o\`u $\hat \mathcal O_{X,x}$ 
d\'esigne le compl\'et\'e de l'anneau local de $X$ en $x$. 
\end{defn}
\par
Soit $X$ une courbe plate et de type fini sur $R$, de fibre g\'en\'erique 
lisse sur $K$. 
Si $x$ est un point lisse de la fibre sp\'eciale de $X$, 
le compl\'et\'e 
de l'anneau local de $X$ en $x$ est isomorphe \`a la $R$-alg\`ebre $\disc$ 
des s\'eries formelles en une variable \`a coefficients dans $R$. 
Si maintenant $x$ est un point double ordinaire de la fibre sp\'eciale de $X$, 
il existe un entier $e$ strictement positif, qu'on appelle {\bf \'epaisseur} 
du point double $x$, tel que le compl\'et\'e 
de l'anneau local de $X$ en $x$ est isomorphe \`a la $R$-alg\`ebre 
$\kuro$. Ceci nous conduit \`a poser les d\'efinitions suivantes :
\begin{defn}
On appelle {\bf disque formel} sur $R$ le $R$-sch\'ema 
$\mathcal D:=\spec \disc$ et 
{\bf couronne formelle d'\'epaisseur} $e \in \mathbb N_{>0}$ le 
$R$-sch\'ema $\mathcal C_e:=\spec \kuro$.
\end{defn}
\subsection{G\'eom\'etrie du disque formel}
Le disque formel sur $R$ est un $R$-sch\'ema lisse. Sa fibre sp\'eciale est le 
germe analytique d'un point lisse d'une courbe alg\'ebrique sur $k$. 
Consid\'erons le disque $D:=\{z \in K^a|v_K(z)>0\}$, si $f$ appartient \`a 
$\disc$, pour tout point $z$ de $D$, la s\'erie $f(z)$ converge. On obtient 
alors une application $\Phi_{\mathcal D}$ de $D$ dans la fibre 
g\'en\'erique de $\mathcal D$, 
qui \`a $z$ associe l'id\'eal premier form\'e des $f$ dans $\disc$ qui 
s'annulent en $z$.
\begin{lem}
L'application $\Phi_{\mathcal D}$ induit une bijection de $D/G_K$ sur 
la fibre g\'en\'erique de $\mathcal D$.
\end{lem}

\begin{figure}[htpb]
         \input{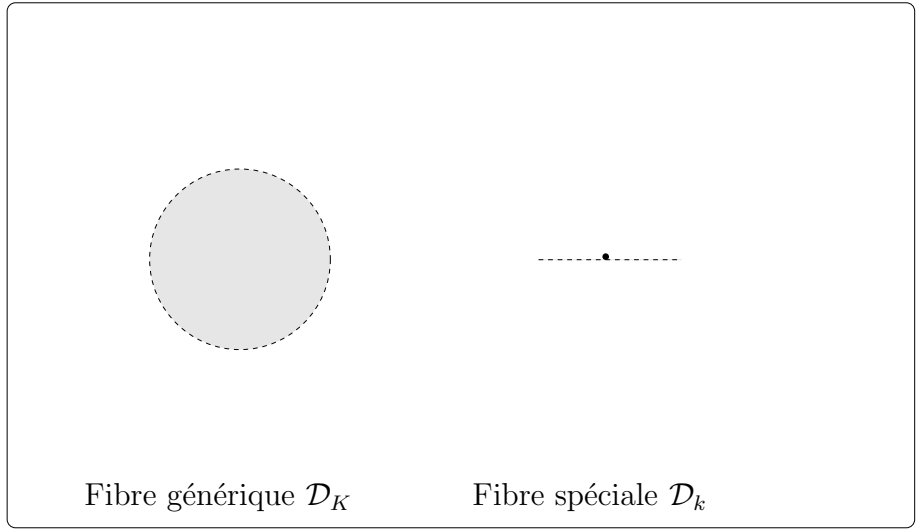}
         \caption{Le disque formel $\mathcal D$ sur $R$}
\end{figure}

Soit $d$ un entier strictement positif. 
Rappelons qu'un polyn\^ome $\displaystyle{\sum_{i=0}^d c_iZ^i}$ unitaire de 
degr\'e $d$ \`a coefficients dans $R$ est dit distingu\'e si
 $\pi$ divise $c_i$ pour $0 \le i<d$.
\begin{proof}
Soit $\frak p$ un id\'eal premier de $\disc$ ne contenant pas $\pi$, le 
th\'eor\`eme de pr\'eparation de Weierstrass entra\^\i ne que $\frak p$ est 
principal, engendr\'e par un polyn\^ome distingu\'e irr\'eductible $P$. 
Si $z$ est une racine de $P$ dans $K^a$, $z$ appartient \`a $D$ car $P$ est 
distingu\'e, et $\Phi_{\mathcal D}(z)=\frak p$. Ainsi $\Phi_{\mathcal D}$ 
est surjective. Par ailleurs, il est clair que $\Phi_{\mathcal D}$ passe au 
quotient. De plus, si $z'$ appartient \`a $D$ et 
$\Phi_{\mathcal D}(z')=\frak p$, alors $z'$ est une racine de $P$ et donc 
$z'$ est un conjugu\'e de $z$ sous l'action de $G_K$.
\end{proof}
\newpage
\subsection{G\'eom\'etrie de la couronne formelle d'\'epaisseur $e$}\ \\
\par
L'anneau $\kouro:=\kuro$ est local complet, d'id\'eal maximal 
$(\pi, Z_1, Z_2)$, 
int\`egre, noeth\'erien et de dimension $2$. Un \'el\'ement $f$ de $\kouro$ 
s'\'ecrit de mani\`ere unique comme une s\'erie de Laurent \`a coefficients 
dans $R$
$$(*)\qquad f:=\sum_{\nu \geq 0} a_{\nu} Z_1^\nu + 
\sum_{\nu >0} a_{-\nu}Z_2^{\nu}.$$
Cette remarque conduit \`a la
\begin{defn}
On appelle {\bf coordonn\'ee de Laurent} sur la couronne formelle 
d'\'epaisseur $e$ un \'el\'ement $Z$ de $\kouro$ tel que 
$\displaystyle{\frac {\pi^e}{Z}}$ appartienne \`a $\kouro$ et, 
pour tout $f$ dans $\kouro$, il existe une unique famille 
$(a_{\nu})_{\nu \in \mathbb Z}$ telle 
que
$$f:=a_0+\sum_{\nu>0}(a_{\nu} Z^\nu + a_{-\nu}(\frac{\pi^e}{Z})^{\nu}).$$
Autrement dit, si $Z':=\displaystyle{\frac {\pi^e}{Z}}$, on a 
$\kouro=\displaystyle{\frac {R[[Z,Z']]}{(ZZ'-\pi^e)}}$.
\end{defn}
\begin{defn}
\par
On appellera {\bf bord} de la couronne formelle le point g\'en\'erique 
d'une composante irr\'eductible de la fibre sp\'eciale de $\mathcal C_e$. 
Une couronne formelle poss\`ede donc deux bords. 
Si $Z$ est une coordonn\'ee de Laurent sur $\mathcal C_e$, alors l'id\'eal 
premier $\frak p_Z:=(\pi,\displaystyle{\frac {\pi^e}{Z}})$ de $\kouro$ est 
un bord de $\mathcal C_e$. On dira que $\frak p_Z$ est le bord 
correspondant \`a $Z$.
\par
Si $\eta$ est un bord, l'anneau local $\mathcal O_{\mathcal C_e,\eta}$ 
est un anneau de valuation discr\`ete, \text{d'uniformisante} $\pi$. On notera 
$v_{\eta}$ la valuation correspondante du corps des fractions $\mathcal K_e$ 
de $\kouro$ qui prolonge $v_K$. Le corps r\'esiduel de 
$(\mathcal K_e,v_{\eta})$ est un corps de s\'eries de Laurent en une variable 
sur $k$, qu'on notera $k((\eta))$. Si $Z$ est une coordonn\'ee de Laurent avec 
$\eta=\frak p_Z$, on notera \'egalement $v_Z=v_{\eta}$ ; on a alors 
$k((\eta))=k((z))$, o\`u $z=Z \mod \pi$. 
On v\'erifie que si $f$ appartient \`a $\mathcal A_e$, 
\text{$f=a_0+\sum_{\nu>0}(a_{\nu} Z^\nu + a_{-\nu}(\frac{\pi^e}{Z})^{\nu})$,} 
alors
$$v_{Z}(f)=\min_{\nu \in \mathbb{Z}}(v_K(a_\nu)+e\max(0;-\nu)).$$
On remarquera que le compl\'et\'e de $\mathcal O_{\mathcal C_e,\eta}$ 
s'identifie \`a $R[[Z]]\{Z^{-1}\}$.
\par
Le corps r\'esiduel $k((\eta))$ de $(\mathcal K_e,v_{\eta})$ est muni 
d'une valuation discr\`ete $\ord_{\eta}$, normalis\'ee de fa\c con \`a ce 
qu'une uniformisante de $k((\eta))$ ait pour valuation $1$. 
Si $f$ est non nul dans $\mathcal K_e$, 
$f$ s'\'ecrit de mani\`ere unique $f=\pi^{v_{\eta}(f)}f_0$, avec 
$f_0$ dans $\mathcal O_{\mathcal C_e,\eta}$. On note alors 
$\ord_{\eta}(f):=\ord_{\eta}(\bar f_0)$, o\`u $\bar f_0$ est 
l'image r\'esiduelle de $f_0$ dans $k((\eta))$. Si $Z$ est une coordonn\'ee 
de Laurent correspondant au bord $\eta$, on v\'erifie que pour 
\text{$f=a_0+\sum_{\nu>0}(a_{\nu} Z^\nu + a_{-\nu}(\frac{\pi^e}{Z})^{\nu})$,}
$$\ord_{Z}(f):=\ord_{\eta}(f)=\min\{\nu \in \mathbb{Z}|v_K(a_\nu)+e\max(0;-\nu)
=v_{Z}(f)\}.$$
L'application $w_Z:\mathcal K_e \to \mathbb Z 
\times \mathbb Z \cup \{\infty\}$, qui \`a $f$ non nul 
associe $(v_\eta(f),ord_\eta(f))$ est une 
valuation de rang $2$, de corps r\'esiduel $k$, en munissant 
$\mathbb Z \times \mathbb Z$ de l'ordre lexicographique.
\end{defn}
 Comme pour le disque formel, le point clef pour comprendre la g\'eom\'etrie 
de la fibre g\'en\'erique d'une couronne formelle est une version ad\'equate 
du th\'eor\`eme de pr\'eparation de Weierstrass, que nous donnons ci-dessous :
\begin{lem}
\label{thprep}
Soient $Z$ une coordonn\'ee de Laurent de $\mathcal C_e$, 
$Z'=\frac{\pi^e}{Z}$, et 
$f \in \kouro$, non inversible, 
avec $v_{Z'}(f)=0$, (en particulier, $\ord_{Z'}(f)>0$) ; 
il existe alors un polyn\^ome distingu\'e $P$ \`a coefficients dans $R$ et 
un inversible $U$ de $\kouro$, tels que : 
$$Z^{\ord_{Z'}(f)}f=\pi^{v_{Z}(f)}P(Z)U.$$
De plus, le degr\'e de $P$ est alors $\ord_{Z}(f)+\ord_{Z'}(f)$.
\end{lem}
Nous donnons une preuve de ce r\'esultat, par manque de r\'ef\'erence 
ad\'equate dans la litt\'erature.

\begin{proof}
Posons $\nu_0:=\ord_{Z'}(f)$, on d\'efinit les endomorphismes de $\kouro$ 
$R$-lin\'eaires $\omega$, $\varphi$, $pr^{>0}$, $pr^{\leq 0}$ par :
\begin{itemize}
\item
$\omega(\sum_{\nu \geq 0}b_{\nu}Z^{\nu}+\sum_{\nu >0}b_{-\nu}{Z'}^{\nu})
:=\sum_{\nu \geq 0}b_{-(\nu +\nu_0)}{Z'}^{\nu}$
\item
$\varphi(\sum_{\nu \geq 0}b_{\nu}Z^{\nu}+\sum_{\nu >0}b_{-\nu}{Z'}^{\nu})
:=\sum_{\nu \geq 0}b_{\nu}Z^{\nu}+\sum_{0<\nu<\nu_0}b_{-\nu}{Z'}^{\nu}$
\item
$pr^{>0}(\sum_{\nu \geq 0}b_{\nu}Z^{\nu}+\sum_{\nu >0}b_{-\nu}{Z'}^{\nu})
:=\sum_{\nu>0}b_{\nu}Z^{\nu}$
\item
$pr^{\leq 0}:=id_{\kouro}-pr^{>0}$
\end{itemize}
Remarquons que, pour tout  $g \in \kouro$, $g=\varphi(g)+{Z'}^{\nu_0}\omega(g)$. 
Montrons qu'il existe $q \in \kouro$ tel que $\omega(qf)=1$. Pour $f$ dans 
$\kouro$, $\omega(qf)=\omega(q\varphi(f))+\omega({Z'}^{\nu_0}
q\omega(f))$. 
Or, si $h \in \kouro,\ \omega({Z'}^{\nu_0}h)=pr^{\leq 0}(h)$. Il suit :
$$\forall q \in \kouro\ \ \omega(qf)=q\omega(f)+\omega(q \varphi(f))
-pr^{>0}(q\omega(f)).$$
Soit $\theta :\kouro \rightarrow \kouro$ l'op\'erateur $\omega \circ \frac 
{\varphi(f)}{\omega(f)}\ -pr^{>0}$, l'\'el\'ement $\omega(f)$ \'etant 
inversible dans $\kouro$, trouver $q$ tel que $\omega(qf)=1$ revient
\`a trouver $W:=q\omega(f)$ tel que $(id_{\kouro} +\theta)(W)=1$. Par 
d\'efinition de $\nu_0$, $\varphi(f) \in (\pi,Z)$, donc 
$\frac {\varphi(f)}{\omega(f)} \kouro \subset (\pi,Z)$. 
Posons $\delta:=\omega \circ \frac {\varphi(f)}{\omega(f)}$. On a $\theta(1)
=\delta(1)$. Comme $pr^{>0} \circ \omega=0$, il suit par r\'ecurrence 
sur $n$ que pour tout entier positif $n$, $\theta^n(1)=\delta^n(1)$. 
Ecrivons $\frac {\varphi(f)}{\omega(f)}:=\pi f_1+Zf_2$, $\omega(Zf_2)$ est
 divisible par $\pi^e$ dans $\kouro$. Donc 
$\delta(1) \in \pi \kouro$. Par r\'ecurrence sur $n \in \mathbb{N}$, 
$\theta^n(1)=\delta^n(1) \in \pi^n\kouro$. Posons alors 
$W:=\sum_{h \geq 0}(-1)^h \theta^{\circ h}(1)$. On a 
$(id_{\kouro} +\theta)(W)=1$. 
Soit $q:=W\omega(f)^{-1}$. On a alors $\omega(qf)=1$. En regardant le coefficient de ${Z'}^{\nu_0}$ dans $qf$, on en d\'eduit ais\'ement que $q$ est 
inversible. Donc, $f=q^{-1}h$, o\`u les coefficients de ${Z'}^{\nu}$ pour 
$\nu>\nu_0$ dans $h$ sont tous nuls, i.e. $Z^{\nu_0}h \in R[[Z]]$. Le 
th\'eor\`eme de pr\'eparation de Weierstrass dans $R[[Z]]$ permet donc 
d'\'ecrire 
$Z^{\nu_0}h=\pi^rPu_1$, avec $r \in \mathbb{N}$, $u_1$ inversible et 
P un polyn\^ome distingu\'e. Si $U:=q^{-1}u_1$, on a $Z^{\nu_0}f=\pi^rPU$. 
Comme $P$ est unitaire en $Z$, on a $v_{Z}(P)=0$. Donc $v_{Z}(f)=r$. 
De plus, le degr\'e de $P$ vaut :
$$\deg (P)=\ord_{Z}(P)=\ord_{Z}(f)+\nu_0=\ord_{Z}(f)+\ord_{Z'}(f).$$
\end{proof}
Notons $C_e:=\{z \in K^a|0<v_K(z)<e\}$, le choix d'une coordonn\'ee de Laurent $Z$ sur $\mathcal C_e$ permet de d\'efinir une application 
$\Phi_{\mathcal C_e}$ de $C_e$ dans la fibre g\'en\'erique de 
$\mathcal C_e$ par $\Phi_{\mathcal C_e}(z):=\{f \in \kouro|f(z)=0\}$. 
\begin{cor}
L'anneau $\kouro\otimes_R K$ est principal. De plus, l'application 
$\Phi_{\mathcal C_e}$ induit une bijection de $C_e/G_K$ dans la fibre 
g\'en\'erique de $\mathcal C_e$
\end{cor}
\begin{figure}[htpb]
         \input{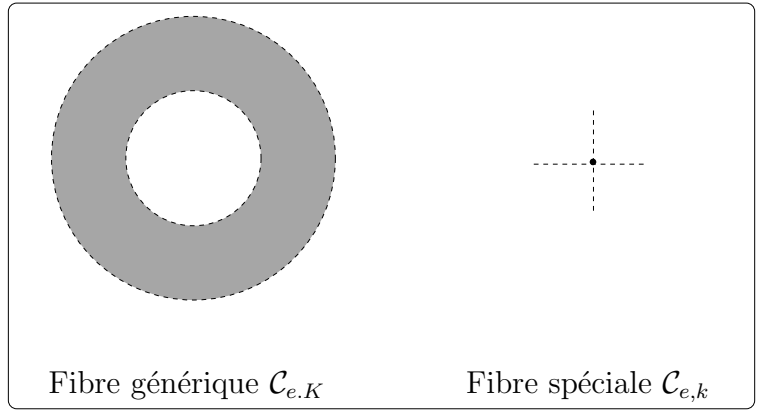}
         \caption{La couronne formelle $\mathcal C_e$ sur $R$}
    \end{figure}

\section{Techniques de recollement formel}
\par
Consid\'erons une $R$-courbe formelle nodale $\mathcal X$, si $x$ est un 
point ferm\'e de la fibre sp\'eciale $\mathcal X_k:=\mathcal X \times_{R}k$, 
$\mathcal X$ peut alors \^etre vue comme recollement de la fibre formelle 
de $x$ et de l'ouvert $\mathcal X\setminus \{x\}$. Les donn\'ees du 
recollement sont inscrites 
dans une suite exacte de $R$-modules (\ref {ser}). Cette remarque est le 
point de d\'epart d'une technique de construction de rel\`evements sur 
$R$ de rev\^etements de courbes nodales sur $k$ (\ref{recolloc}).
\subsection{Un lemme d'alg\`ebre topologique}
 Le lemme suivant, de preuve imm\'ediate, sera utilis\'e constamment 
dans la suite. On fixe une uniformisante $\pi$ de $R$. 
Par un $R$-module $M$ complet pour la topologie $\pi$-adique, nous 
sous-entendrons ici un module s\'epar\'e et complet. Pour 
tout $R$-module $M$, on d\'esigne par $\overline M$ le $k$-espace vectoriel 
$M\otimes_Rk$.
\begin{lem}
\label{lat}
(i) Soit $M_1 \stackrel{\phi}{\longrightarrow} M_2$ un homomorphisme de 
$R$-modules, on suppose $M_1$ complet et $M_2$ s\'epar\'e pour la topologie 
$\pi$-adique. 
Si l'homomorphisme de $k$-espaces vectoriels $\overline{M_1} \stackrel{\overline{\phi}}{\longrightarrow} \overline{M_2}$ est surjectif, alors $\phi$ est surjectif.\\
(ii) Sous les m\^emes hypoth\`eses, si de plus $M_2$ est plat sur $R$ et $\overline{\phi}$ est un isomorphisme de $k$-espaces vectoriels, alors $\phi$ est un isomorphisme de $R$-modules.\\
(iii) Soient $u:M_1 \longrightarrow M_2$ et 
$v:M_2\longrightarrow M_3$ des homomorphismes de $R$-modules tels que $v \circ u=0$, on suppose $M_1$ et $M_2$ complets, et $M_3$ s\'epar\'e pour la topologie $\pi$-adique. On demande \'egalement que les modules $M_2$ et $M_3$ soient plats sur $R$, c'est-\`a-dire sans torsion. 
Si la suite de $k$-espaces vectoriels $0 \longrightarrow \overline{M_1} \stackrel{\overline{u}}{\longrightarrow} \overline{M_2} 
\stackrel{\overline{v}}{\longrightarrow} \overline{M_3} \longrightarrow 0$
est exacte, il en est de m\^eme de la suite de $R$-modules 
$$0 \longrightarrow M_1 
\stackrel{u}{\longrightarrow} M_2  \stackrel{v}{\longrightarrow}
 M_3 \longrightarrow 0.$$
\end{lem}
\subsection{Une suite exacte de recollement}
Pour $\mathcal X$ un $R$-sch\'ema formel, on notera 
$\irr(\mathcal X)$ l'ensemble des points g\'en\'eriques des 
composantes irr\'eductibles de $\mathcal X$. Pour simplifier, si $A$ est 
une $R$-alg\`ebre compl\`ete pour la topologie $\pi$-adique, on notera 
$\irr(A):=\irr(\spf A)$. 
Soit $\mathcal X$ un $R$-sch\'ema formel, de dimension $1$, 
localement noeth\'erien, plat sur $R$.  Soient $x_1,\dots,x_r$ 
des points ferm\'es de $\mathcal X$, et 
$\mathcal U=\mathcal X \setminus \{x_1,\dots,x_r\}$. Le 
diagramme de $\mathcal O_{\mathcal X}$-modules 
ci-dessous est commutatif : 
$$
  \diagram
       \mathcal{O}_{\mathcal U}\drto^{l'_1}
    & \mathcal{O}_{\mathcal X} \lto_{\rho} \rto^{l_0} \dto^{l_1}
    & \displaystyle{\prod_{q=1}^r} \hat{\mathcal{O}}_{\mathcal{X},x_q}
       \dto^{l_2}\\
    & \displaystyle{\prod_{\xi \in \irr \mathcal{X}}}
     \hat{\mathcal{O}}_{\mathcal{X},\xi}    \rto^{l_3 \qquad} 
    & \displaystyle{\prod_{q=1}^r \prod_{\xi \in \irr(\hat{\mathcal O}_{\mathcal X,x_q})}}
       (\hat{\mathcal O}_{\mathcal X,x_q})_{\xi}^{\wedge}
  \enddiagram 
$$

\begin{prop}
\label{ser}
Soient $\mathcal{X}$ un sch\'ema formel localement noeth\'erien, 
de dimension $1$, plat sur $R$, $x_1,\dots,x_r$ 
des points ferm\'es de $\mathcal X$ et 
$\mathcal U:=\mathcal X \setminus \{x_1,\dots,x_r\}$. On a alors la 
suite exacte de $\mathcal O_{\mathcal X}$-modules :
$$
0 \to \mathcal O_{\mathcal X} \stackrel{\Delta}{\to} \mathcal O_{\mathcal U} 
\times \prod_{q=1}^r \hat{\mathcal{O}}_{\mathcal{X},x_q} \stackrel{\theta}{\to}
    \prod_{q=1}^r \prod_{\xi \in \irr(\hat{\mathcal O}_{\mathcal X,x_q})}
      (\hat{\mathcal O}_{\mathcal X,x_q})_{\xi}^{\wedge} \to 0
$$
o\`u $\Delta=(\rho,l_0)$ est le morphisme diagonal et 
$\theta(g,\{f_q\})=(l_3 \circ l_1'(g)-l_2(\{f_q\}))$.
\end{prop}
\begin{rem}
Cette suite exacte apparait dans les travaux d'Harbater et 
Stevenson (\cite{H-S} section 1) dans le contexte d'\'egale 
caract\'eristique.
\end{rem}
\begin{proof}On a $\theta \circ \Delta=0$ par la commutativit\'e du diagramme 
ci-dessus. 
Notons $X$ (resp. $U$) la fibre sp\'eciale de $\mathcal X$ 
(resp. $\mathcal U$). 
D'apr\`es le lemme \ref{lat}, il suffit de montrer que la suite ci-dessous 
obtenue en tensorisant par $k$ est exacte :
$$
0 \to \mathcal O_{X}\to \mathcal O_U \times \prod_{q=1}^r 
\hat{\mathcal O}_{X,x_q} \to \prod_{q=1}^r 
\prod_{\xi \in \irr(\hat{\mathcal O}_{X,x_q})}
(\hat{\mathcal O}_{X,x_q})_{\xi}^{\wedge} \to 0
$$
Comme $X$ est localement noeth\'erien, il suffit de tester l'exactitude 
au niveau des fibres compl\'et\'ees. La v\'erification est alors imm\'ediate.
\end{proof}
\begin{exmp}
Si $\mathcal X$ est le disque unit\'e ferm\'e standard, c'est-\`a-dire 
$\mathcal X=\spf R\{T\}$, et $x:=(\pi,T)$, on obtient la suite exacte
$$0\longrightarrow R\{T\}\stackrel{\Delta}{\longrightarrow} R\{T,T^{-1}\}\times R[[T]]\stackrel{-}\longrightarrow R[[T]]\{T^{-1}\}\longrightarrow0.$$
\end{exmp}
\begin{exmp}
Si $\mathcal X$ est la couronne ferm\'ee standard d'\'epaisseur $e$, 
c'est-\`a-dire 
$\mathcal X=\spf \displaystyle{\frac{R\{T_1,T_2\}}{(T_1T_2-\pi^e)}}$, 
et $x:=(\pi,T_1,T_2)$, on obtient la suite exacte 
$$0\rightarrow \frac{R\{T_1,T_2\}}{(T_1T_2-\pi^e)}\stackrel{\Delta}
{\rightarrow} \prod_{i=1,2}R\{T_i,T_i^{-1}\}\times 
\frac{R[[T_1,T_2]]}{(T_1T_2-\pi^e)}
\stackrel{\theta}{\rightarrow}
 \prod_{i=1,2}R[[T_i]]\{T_i^{-1}\} \rightarrow 0$$
o\`u $\theta(f_1,f_2,g)=(f_1-g,f_2-g)$ (avec les identifications 
\'evidentes en termes de s\'erie de Laurent).
\end{exmp}
\subsection{Construction d'un mod\`ele avec \'epaisseurs fix\'ees}
Soit $X_0$ une courbe nodale projective sur $k$, on appelle 
mod\`ele de $X_0$ sur $R$ un couple $(X,\phi)$, o\`u $X$ est un $R$-sch\'ema 
propre et plat, et $\phi: X \times_R k \to X_0$ est un isomorphisme de 
$k$-sch\'emas. 
\begin{prop}
\label{modele}
Soit $X_0$ une courbe nodale projective sur $k$, $S$ l'ensemble des 
points ferm\'es singuliers de $X_0$, et $(e_x)_{x \in S}$ une famille 
d'entiers strictement positifs. Il existe alors un mod\`ele $(X,\phi)$ de 
$X_0$ sur $R$ tel que l'\'epaisseur de $x$ dans $X$ soit $e_x$, pour tout 
point $x$ de $S$.
\end{prop}
\begin{proof}
Lorsque $X_0$ est lisse sur $k$, la proposition r\'esulte du corollaire 7.4 
de \cite{SGA1}, expos\'e III. Supposons \`a pr\'esent $X_0$ irr\'eductible. 
Soit $\tilde X_0$ la normalis\'ee de $X_0$, elle poss\`ede un mod\`ele 
$(\tilde X,\tilde \phi)$ sur $R$. Notons, pour $x$ dans $S$, 
$x_1$ et $x_2$ les points de $\tilde X_0$ au-dessus de $x$, et 
$U_0$ l'ouvert $X_0 \setminus S$ de $X_0$, qui s'identifie \`a un ouvert de 
$\tilde X_0$. Notons $\mathcal U$ le sous-sch\'ema formel ouvert de 
$\tilde \mathcal X$ (compl\'et\'e de $\tilde X$ le long de la fibre 
sp\'eciale) qui correspond \`a $U_0$. Si 
$\frac{R[[a_x,b_x]]}{(a_xb_x-\pi^{e_x})}$ et 
$(\mathcal O_{\tilde X,x_1})_{(\pi)}^{\wedge}\times
(\mathcal O_{\tilde X,x_2})_{(\pi)}^{\wedge}$ 
sont vus comme faisceau gratte-ciel en $x\in S$, consid\'erons 
l'homomorphisme de faisceau de $R$-modules sur $X_0$:
$$
\theta : \mathcal O_{\mathcal U}\times \prod_{x \in S} 
\frac{R[[a_x,b_x]]}{(a_xb_x-\pi^{e_x})} \to 
\prod_{x \in S}((\mathcal O_{\tilde X,x_1})_{(\pi)}^{\wedge}\times
(\mathcal O_{\tilde X,x_2})_{(\pi)}^{\wedge})
$$
obtenu en identifiant $(\mathcal O_{\tilde X,x_1})_{(\pi)}^{\wedge}$ 
(resp. $(\mathcal O_{\tilde X,x_2})_{(\pi)}^{\wedge}$) \`a 
$R[[a_x]]\{a_x^{-1}\}$ (resp. $R[[b_x]]\{b_x^{-1}\}$). Le lemme \ref{lat} 
montre que $\theta$ est surjectif. Son noyau, 
not\'e $\mathcal O_{\mathcal X}$ est un faisceau de $R$-alg\`ebres 
sur l'espace topologique $X_0$, et on voit que l'espace localement annel\'e 
$(X_0,\mathcal O_{\mathcal X})$ est un $R$-sch\'ema formel propre et plat, 
de fibre 
sp\'eciale projective. C'est donc le compl\'et\'e d'un $R$-sch\'ema propre 
et plat le long de la fibre sp\'eciale, qui fournit un mod\`ele pour $X_0$. 
Si maintenant on ne suppose plus $X_0$ irr\'eductible, soit $S'$ le 
sous-ensemble de $S$ form\'e des $x$ qui appartiennent \`a deux 
composantes irr\'eductibles distinctes et $U'_0$ l'ouvert $X_0 \setminus S'$, 
on choisit un mod\`ele $\tilde X_i$ pour chaque composante irr\'eductible 
$X_{0,i}$ de $X_0$. Comme ci-dessus, on construit un homomorphisme surjectif 
de faisceaux de $R$-modules sur $X_0$
$$
\prod_{i}\mathcal O_{\tilde \mathcal X_i}\times \prod_{x \in S} 
\frac{R[[a_x,b_x]]}{(a_xb_x-\pi^{e_x})} \to 
\prod_{x \in S}((\mathcal O_{\tilde X,x_1})_{(\pi)}^{\wedge}\times
(\mathcal O_{\tilde X,x_2})_{(\pi)}^{\wedge})
$$
dont le noyau $\mathcal O_{\mathcal X}$ est un faisceau de $R$-alg\`ebres 
sur l'espace topologique $X_0$. Le $R$-sch\'ema formel propre et plat 
$(X_0,\mathcal O_{\mathcal X})$ s'alg\'ebrise comme ci-dessus, et on 
obtient ainsi le mod\`ele d\'esir\'e.
\end{proof}
\begin{rem}
Le lecteur trouvera une preuve rigide dans \cite{Sa1}, 
lemme 6.3.
\end{rem}
\subsection[Principe local-global]{Principe local-global pour les 
rev\^etements de courbes formelles}
Soit $\overline f:Y \longrightarrow X$ un morphisme s\'eparable fini
entres courbes alg\'ebriques sur $k$, connexes, affines,
r\'eduites. Soient $x$ un point ferm\'e de $X$ et $X'$ l'ouvert
compl\'ementaire du point $x$. On suppose $X'$ lisse sur $k$, 
$\overline f$ \'etale au-dessus de $X'$ et l'image r\'eciproque de $x$ 
r\'eduite \`a un point ferm\'e $y$. On se placera dans l'une des deux 
situations suivantes :\\
 (A) $x$ (resp. $y$) est un point lisse de $X$ (resp. $Y$).
\par
Soient $t$ (resp. $z$) une uniformisante de $X$ (resp. $Y$) en
$x$ (resp. $y$). On suppose
$$\hat{\mathcal{O}}_{Y,y}=k[[z]]=k[[t]][z]/(\overline{P}(z))$$
 o\`u $\overline{P}$ est un polyn\^ome d'Eisenstein s\'eparable de $k[[t]]$.\\
 (B) $x$ (resp. $y$) est un point double ordinaire de $X$
(resp. $Y$),\\
$$\hat{\mathcal{O}}_{X,x}=\frac {k[[t_1,t_2]]}{(t_1t_2)} \text{ et
}\hat{\mathcal{O}}_{Y,y}=\frac{k[[z_1,z_2]]}{(z_1z_2)}=\frac{\hat{\mathcal{O}}_{X,x}[z_1,z_2]}{(\overline{P}(z_1),\overline{Q}(z_2),z_1z_2)},$$
o\`u $\overline{P}$ (resp. $\overline{Q}$) est un polyn\^ome
d'Eisenstein s\'eparable de $k[[t_1]]$ (resp. $k[[t_2]]$).

\begin{thm} {\bf (Principe local-global formel)} 
\label{recolloc}
Soit $\mathcal{X}$ un sch\'ema formel affine normal, plat et topologiquement de type fini sur $R$, de fibre sp\'eciale $X$. On note $\mathcal{X'}$ l'ouvert de $\mathcal{X}$ correspondant \`a $X'$. La restriction de $\overline f$ au-dessus de $X'$ s'\'etend de mani\`ere unique (\`a isomorphisme pr\`es) en un rev\^etement \'etale $f':\mathcal{Y'} \longrightarrow \mathcal{X'}$.
\par
On se donne une $\hat{\mathcal{O}}_{\mathcal{X},x}$-alg\`ebre $A$ finie, normale, $R$-plate et un diagramme commutatif \`a lignes exactes :
$$
  \diagram
    0 \rto &\pi A \rto &A \rto &\hat{\mathcal O}_{Y,y}\rto & 0\\
    0 \rto &\pi \hat{\mathcal O}_{\mathcal X,x} \uto \rto &
    \hat{\mathcal O}_{\mathcal X,x} \rto \uto &\hat{\mathcal O}_{X,x}\uto 
    \rto & 0
  \enddiagram  
$$
 (i) Il existe alors un rev\^etement fini $f:\mathcal{Y}\longrightarrow \mathcal{X}$ relevant $\overline f$ tel que $\mathcal{Y}$ est normal, 
$f_{|\mathcal{X'}}=f'$ et $f$ induit l'extension $\hat{\mathcal{O}}_{\mathcal{X},x}\longrightarrow A$.\\
 (ii) Si de plus $\overline f$ est galoisien de groupe de Galois $G$, et si 
$A$ est munie d'une action de $G$ de sorte que 
$A^G=\hat{\mathcal{O}}_{\mathcal{X},x}$ et que l'homomorphisme de 
$R$-alg\`ebres $A \longrightarrow \hat{\mathcal{O}}_{Y,y}$ soit $G$-\'equivariant,
 le rev\^etement $f:\mathcal{Y}\longrightarrow \mathcal{X}$ est 
galoisien, de groupe de Galois $G$, relevant l'action sur $Y$.
\end{thm}
\begin{proof} On commence par traiter le cas (A) :\\
(i) La $R$-alg\`ebre $(\hat{\mathcal O}_{\mathcal X,x})_{(\pi)}^\wedge$ est isomorphe 
\`a $R[[T]]\{T^{-1}\}$. C'est donc un anneau de valuation 
discr\`ete complet. Comme $A\otimes_{\hat{\mathcal O}_{\mathcal X,x}}
(\hat{\mathcal O}_{\mathcal X,x})_{(\pi)}^\wedge$ est fini sur $(\hat{\mathcal O}_{\mathcal X,x})_{(\pi)}^\wedge$, il est semi-local complet. On a de plus 
le diagramme commutatif :
$$
  \diagram 
    A\otimes_{\hat{\mathcal O}_{\mathcal X,x}}(\hat{\mathcal O}
      _{\mathcal X,x})_{(\pi)}^\wedge  \rto^{\mod \pi\quad}
      &Fr(\hat{\mathcal O}_{Y,y})=k((z))\\
    (\hat{\mathcal O}_{\mathcal X,x})_{(\pi)}^\wedge\rto^{\mod \pi\quad}\uto 
      &Fr(\hat{\mathcal O}_{X,x})=k((t))\uto
  \enddiagram  
$$
Il suit que $A\otimes_{\hat{\mathcal O}_{\mathcal X,x}}(\hat{\mathcal O}
_{\mathcal X,x})_{(\pi)}^\wedge$ est local complet, noeth\'erien, 
d'id\'eal maximal 
engendr\'e par $\pi$. C'est donc un anneau de valuation discr\`ete complet. 
Soit $P$ un polyn\^ome unitaire de $(\hat{\mathcal O}_{\mathcal X,x})_{(\pi)}^\wedge[X]$ relevant $\overline{P}\in k((t))[X]$. L'anneau local 
$A\otimes_{\hat{\mathcal O}_{\mathcal X,x}}(\hat{\mathcal O}_{\mathcal X,x})_{(\pi)}^\wedge$ \'etant hens\'elien, la racine $z$ de $\overline{P}$ 
dans $k((z))$ se rel\`eve en une racine $Z$ de $P$ dans $A\otimes_{\hat{\mathcal O}_{\mathcal X,x}}(\hat{\mathcal O}_{\mathcal X,x})_{(\pi)}^\wedge$. L'homomorphisme de $(\hat{\mathcal O}_{\mathcal X,x})_{(\pi)}^\wedge$-alg\`ebres 
$$\frac{(\hat{\mathcal O}_{\mathcal X,x})_{(\pi)}^\wedge[X]}{(P(X))} \stackrel{u}{\longrightarrow} A\otimes_{\hat{\mathcal O}_{\mathcal X,x}}(\hat{\mathcal O}_{\mathcal X,x})_{(\pi)}^\wedge$$ qui envoie $X$ sur $Z$ est un 
isomorphisme : Il nous suffit d'appliquer le lemme \ref{lat}(ii), car c'est vrai au niveau r\'esiduel par hypoth\`ese.
\par
L'anneau $\mathcal{O}(\mathcal Y') \otimes_{\mathcal{O}(\mathcal X')}
(\hat{\mathcal O}_{\mathcal X,x})_{(\pi)}^\wedge$ est fini sur $(\hat{\mathcal O}_{\mathcal X,x})_{(\pi)}^\wedge$, donc semi-local complet.
\begin{eqnarray*}
\frac{\mathcal{O}(\mathcal Y') \otimes_{\mathcal{O}(\mathcal X')}
(\hat{\mathcal O}_{\mathcal X,x})_{(\pi)}^\wedge}{(\pi)} 
\simeq & &\overline{\mathcal{O}(\mathcal Y')} \otimes_{\mathcal{O}(\mathcal X')}(\hat{\mathcal O}_{\mathcal X,x})_{(\pi)}^\wedge\\
=& & \overline{\mathcal{O}(\mathcal Y')} \otimes_
{\overline{\mathcal{O}(\mathcal X')}}\overline{\mathcal{O}(\mathcal X')}
\otimes_{\mathcal{O}(\mathcal X')}(\hat{\mathcal O}_{\mathcal X,x})_{(\pi)}^\wedge\\
=& & \mathcal{O}(Y')\otimes_{\mathcal{O}(X')}k((t))\\
=& & k((z))
\end{eqnarray*}
 L'anneau $\mathcal{O}(\mathcal Y') \otimes_{\mathcal{O}(\mathcal X')}
(\hat{\mathcal O}_{\mathcal X,x})_{(\pi)}^\wedge$ est donc en fait 
local complet, noeth\'erien, d'id\'eal maximal engendr\'e par $\pi$. C'est donc un anneau de valuation discr\`ete complet. Comme ci-dessus, on en d\'eduit l'isomorphisme
 de $(\hat{\mathcal O}_{\mathcal X,x})_{(\pi)}^\wedge$-alg\`ebres 
$$\frac{(\hat{\mathcal O}_{\mathcal X,x})_{(\pi)}^\wedge[X]}{(P(X))} 
\stackrel{w}{\longrightarrow}\mathcal{O}(\mathcal Y') \otimes_{\mathcal{O}(\mathcal X')}(\hat{\mathcal O}_{\mathcal X,x})_{(\pi)}^\wedge$$ qui envoie $X$ sur $Z'$, o\`u $Z'$ est une racine de $P$ dans $\mathcal{O}(\mathcal Y') \otimes_{\mathcal{O}(\mathcal X')}(\hat{\mathcal O}_{\mathcal X,x})_{(\pi)}^\wedge$ 
relevant $z$.
\par
Soit $\mu:=u\circ w^{-1}$. C'est un isomorphisme de $(\hat{\mathcal O}_{\mathcal X,x})_{(\pi)}^\wedge$-alg\`ebres :
$$\mathcal{O}(\mathcal Y') \otimes_{\mathcal{O}(\mathcal X')}
(\hat{\mathcal O}_{\mathcal X,x})_{(\pi)}^\wedge \stackrel{\mu}{\simeq}
A\otimes_{\hat{\mathcal{O}}_{\mathcal X,x}}
(\hat{\mathcal O}_{\mathcal X,x})_{(\pi)}^\wedge.$$
Notons $\theta$ l'homomorphisme de $\mathcal O(\mathcal X)$-modules d\'efini ci-dessous :
\begin{eqnarray*}
\mathcal O(\mathcal Y')\times A & & \stackrel{\theta}{\longrightarrow} 
A\otimes_{\hat{\mathcal{O}}_{\mathcal X,x}}
(\hat{\mathcal O}_{\mathcal X,x})_{(\pi)}^\wedge\\
(h',g)& & \mapsto  \mu(h'\otimes 1)-g\otimes 1
\end{eqnarray*}
Remarquons que $\theta$ est surjectif d'apr\`es le lemme \ref{lat}(i). On voit ais\'ement que son noyau $\mathcal A$ est une sous-$\mathcal O(\mathcal X)$-alg\`ebre de $\mathcal O(\mathcal Y')\times A$. Comme $\overline{\mathcal A}=
\mathcal O(Y)$ est fini sur $\overline{\mathcal O(\mathcal X)}=\mathcal O(X)$, le lemme (i) dit que $\mathcal A$ est fini sur $\mathcal O(\mathcal X)$.  En particulier, $\mathcal A$ est topologiquement de type fini sur $R$. En vertu du lemme \ref{lat}(iii), les homomorphismes canoniques ci-dessous sont des isomorphismes.
$$\mathcal A\otimes_{\mathcal O(\mathcal X)} \mathcal O(\mathcal X')\simeq
 \mathcal O(\mathcal Y') \text{ et }\mathcal A\otimes_{\mathcal O(\mathcal X)}
 \hat{\mathcal O}_{\mathcal X,x}\simeq A$$
 Le rev\^etement $f: \mathcal Y:=\spf \mathcal A\longrightarrow \mathcal X$ 
convient.\\
(ii) Pour le cas galoisien, remarquons tout d'abord que $\mathcal X'$ est muni canoniquement d'une action de $G$ relevant l'action sur $X'$. On en d\'eduit une action de $G$ sur $\mathcal{O}(\mathcal Y') \otimes_{\mathcal{O}(\mathcal X')}(\hat{\mathcal O}_{\mathcal X,x})_{(\pi)}^\wedge$. Il suffit de  voir que l'isomorphisme $\mu$ donn\'e ci-dessus est \'equivariant. Comme
 $\overline{\mu}$ est \'egal \`a l'identit\'e de $k((z))$, il est 
$G$-\'equivariant. Ainsi, si $\sigma \in G$, $\mu (\sigma(Z'))=
\sigma(\mu(Z'))\mod \pi$. Comme tous deux sont des racines de $P$, ils sont
 \'egaux, ce qui ach\`eve la preuve de (ii).
\newline
\par Pla\c cons nous \`a pr\'esent dans le cas (B) : 
On remarque que, pour $i=1,2$, $(\hat{\mathcal O}_{\mathcal X,x})_{\eta_i}^\wedge$ est isomorphe \`a l'alg\`ebre $R[[T_i]]\{T_i^{-1}\}$. C'est donc un anneau 
de valuation discr\`ete complet. L'anneau $A\otimes_{\hat{\mathcal O}_{\mathcal X,x}}(\hat{\mathcal O}_{\mathcal X,x})_{\eta_i}^\wedge$ est fini sur 
$(\hat{\mathcal O}_{\mathcal X,x})_{\eta_i}^\wedge$, donc semi-local complet.
 On a de plus les diagrammes commutatifs pour $i=1,2$ :
$$
  \diagram
    A\otimes_{\hat{\mathcal O}_{\mathcal X,x}}(\hat{\mathcal O}_{\mathcal X,x})      _{\eta_i}^\wedge \rto^{\qquad \mod \pi} & k((z_i))\\
    (\hat{\mathcal O}_{\mathcal X,x})_{\eta_i}^\wedge \uto \rto^{\mod \pi}
    & k((t_i)) \uto\\
  \enddiagram
$$
\par
On en d\'eduit comme dans le cas A un isomorphisme de 
$(\hat{\mathcal O}_{\mathcal X,x})_{\eta_i}^\wedge$-alg\`ebres pour $i=1,2$ :
$$\mathcal O(\mathcal Y') \otimes_{\mathcal O(X')}(\hat{\mathcal O}_{\mathcal X,x})_{\eta_i}^\wedge\stackrel{\mu_i}{\simeq} A\otimes_{\hat{\mathcal O}_{\mathcal X,x}}(\hat{\mathcal O}_{\mathcal X,x})_{\eta_i}^\wedge$$
Notons $\theta$ l'homomorphisme $\mathcal O(\mathcal Y')\times A
 \stackrel{\theta}{\longrightarrow} A\otimes_{\hat{\mathcal O}_{\mathcal X,x}}(\hat{\mathcal O}_{\mathcal X,x})_{\eta_1}^\wedge \times A\otimes_{\hat{\mathcal O}_{\mathcal X,x}}(\hat{\mathcal O}_{\mathcal X,x})_{\eta_2}^\wedge$ de $A$-modules d\'efini par :
$$\theta((h',g))=(\mu_1(f'\otimes 1)-g\otimes 1,\mu_2(h'\otimes 1)-
g\otimes 1)$$
 On finit comme pour le cas A en consid\'erant $\mathcal A:=\ker \theta$.
 Le rev\^etement $f: \mathcal Y:=\spf \mathcal A\longrightarrow \mathcal X$ convient.
 Le cas galoisien se traite comme dans le cas (A), en montrant que $\mu_1$ et 
$\mu_2$ sont automatiquement équivariants.
\end{proof}
\begin{rem}
Le cas (A) est trait\'e par des m\'ethodes rigides par B.Green et 
M. Matignon dans \cite{G-M 1} (III 1.1). Toujours dans ce cas, une preuve 
cohomologique a \'et\'e donn\'ee par J. Bertin et A. M\'ezard (\cite{B-M}).
\end{rem}
\begin{prop}
\label{invol}
Soient $\sigma$ une involution de $Y$, 
$X=Y/\langle \sigma \rangle$, $y$ un point double
 de $Y$ fixe par $\sigma$, au-dessus du point $x$. On suppose que $\sigma$ 
permute les 
deux points g\'en\'eriques de $\spec \hat{\mathcal O}_{Y,y}$, $x$ est alors 
r\'egulier. De plus, on peut trouver $z_1,\ z_2,\ t=z_1+z_2$ tels que 
$\hat{\mathcal{O}}_{X,x}=k[[t]],\ \hat{\mathcal{O}}_{Y,y}=
\displaystyle{\frac{k[[z_1,z_2]]}{(z_1z_2)}}$ et $\sigma(z_1)=z_2$.\\
Soit $\mathcal{X}$ un sch\'ema formel affine normal, plat et topologiquement de type fini sur $R$, de fibre sp\'eciale $X$. 
Il existe alors un rev\^etement $2$-cyclique $\mathcal{Y}$ de $\mathcal{X}$
 prolongeant $Y \rightarrow X$.
\end{prop}
\begin{proof}
Remarquons tout d'abord que l'on dispose d'un rel\`evement local 
\'equivariant \'evident (avec $e$ un entier strictement positif)
$$
  \diagram
     \frac{R[[Z_1,Z_2]]}{(Z_1Z_2-\pi^e)} \rto^{\mod \pi \quad} & 
         \frac{k[[z_1,z_2]]}{(z_1z_2)}\\
     R[[T]] \uto^{\psi} \rto^{\mod \pi \quad} & k[[t]] \uto^{\overline{\psi}}
  \enddiagram
$$
avec $\psi (T)=Z_1+Z_2$, et $\sigma(Z_1)=Z_2$, $\sigma(Z_2)=Z_1$.
\par
On note $\mathcal{X'}$ l'ouvert de $\mathcal{X}$ correspondant \`a $X'$.
 Comme ci-dessus, la restriction de $\overline f$ au-dessus de $X'$ s'\'etend de mani\`ere unique (\`a isomorphisme pr\`es) en un rev\^etement \'etale $f':\mathcal{Y'} \longrightarrow \mathcal{X'}$, galoisien de groupe $\langle 
\sigma \rangle$.
\par
 L'anneau $\mathcal O(\mathcal Y')\otimes_{\mathcal O(\mathcal X')}
(\hat{\mathcal O}_{\mathcal X,x})_{(\pi)}^\wedge$ est fini sur
$(\hat{\mathcal O}_{\mathcal X,x})_{(\pi)}^\wedge$, donc semi-local complet. 
Comme de plus $\overline{\mathcal O(\mathcal Y')\otimes_{\mathcal 
O(\mathcal X')}(\hat{\mathcal O}_{\mathcal X,x})_{(\pi)}^\wedge}=\mathcal O
(Y')\otimes_{\mathcal O(X')}k((t))=k((z_1))\times k((z_2))$, $\mathcal O(\mathcal Y')\otimes_{\mathcal O(\mathcal X')}
(\hat{\mathcal O}_{\mathcal X,x})_{(\pi)}^\wedge$ est le produit 
$A_1\times A_2$ de deux anneaux locaux. En fait, $A_i$ est une $R$-alg\`ebre
 qui est un anneau de valuation discr\`ete complet, d'uniformisante $\pi$,
 de corps r\'esiduel $k((z_i))$, pour $i=1,2$. Il suit que $A_i$ est 
isomorphe \`a $R[[Z_i]]\{Z_i^{-1}\}$. On notera $\mu$ l'isomorphisme 
\'equivariant de $\mathcal O(\mathcal Y')\otimes_{\mathcal O(\mathcal X')}
(\hat{\mathcal O}_{\mathcal X,x})_{(\pi)}^\wedge$ sur $R[[Z_1]]\{Z_1^{-1}\}
\times R[[Z_2]]\{Z_2^{-1}\}$ qui envoie un rel\`evement $Z_1'$ de 
$z_1$ sur $Z_1$ et $\sigma(Z_1')$ (qui est un rel\`evement de 
$z_2$) sur $Z_2$ . Notons $\theta$ l'homomorphisme 
$$\mathcal O(\mathcal Y')\times \frac{R[[Z_1,Z_2]]}{(Z_1Z_2-\pi^e)}
 \stackrel{\theta}{\longrightarrow}
R[[Z_1]]\{Z_1^{-1}\}\times R[[Z_2]]\{Z_2^{-1}\}$$
d\'efini par $\theta (g,h)=\mu (g\otimes 1)-(\phi_1(h),\phi_2(h))$, o\`u 
$\phi_i$ d\'esigne l'injection canonique 
$\displaystyle{\frac{R[[Z_1,Z_2]]}{(Z_1Z_2-\pi^e)}}\rightarrow R[[Z_i]]\{Z_i^{-1}\}$. Comme pr\'ec\'edemment, on conclut en consid\'erant la $R$-alg\`ebre 
$\ker \theta$. 
\end{proof}

\section{Construction d'automorphismes de couronnes formelles}
\par
Le paragraphe pr\'ec\'edent montre que pour construire un rel\`evement formel 
sur $R$ d'une courbe nodale sur $k$ avec action d'un groupe fini, 
on doit construire un rel\`evement local pour la fibre formelle d'un 
point ferm\'e 
$y$ fixe sous l'action du groupe. Une telle fibre formelle est un disque 
formel si $y$ est un point lisse et une couronne formelle si $y$ est un 
point double. Le cas des disques formels a \'et\'e 
trait\'e dans \cite{G-M 1}. Ici, nous montrons comment les r\'esultats de 
rel\`evement local contenus dans loc. cit. permettent d'obtenir des 
r\'esultats de rel\`evement local pour des groupes cycliques ou 
di\'edraux agissant sur un point double, en utilisant des m\'ethodes de 
recollement formel.
\subsection{Compl\'ements sur le rel\`evement lisse}
Pour la construction de rel\`evements locaux galoisiens pour un point double, 
nous aurons besoin d'ajuster convenablement des rel\`evements locaux pour 
chacune des deux branches du point double. C'est le r\^ole du lemme qui 
suit :
\begin{lem}
\label{ajustement}
 Notons $K_0$ le corps des fractions de $W(k)$. On suppose l'extension 
$K/K_0$ galoisienne. Soient $r$ un entier 
sup\'erieur \`a $1$ et $\bar \sigma$ un automorphisme 
d'ordre $p^r$ de $k[[z]]$. S'il existe un automorphisme 
d'ordre $p^r$ de $R[[Z]]$ fixant $0$ et relevant $\bar \sigma$, 
alors pour toute racine primitive $p^r$-i\`eme de l'unit\'e $\zeta_{(r)}^l$ 
dans $K^{alg}$, il existe un automorphisme $\sigma$
d'ordre $p^r$ de $R[[Z]]$, fixant $0$ et relevant $\bar \sigma$, 
tel que $\displaystyle{\frac{d\sigma}{dZ}}(0)=\zeta_{(r)}^l$.
\end{lem}
\begin{proof} Le groupe de 
Galois Gal$(K^{alg}/K_0)$ agit sur l'ensemble des automorphismes 
$\sigma$ d'ordre $p^r$ de $R[[Z]]$ fixant $0$ 
par l'action sur les coefficients de la s\'erie $\sigma(Z)$. De plus, 
si $\tau \in$ Gal$(K^{alg}/K_0)$, on a $\bar{\sigma^{\tau}}=\bar{\sigma}$. 
Comme Gal$(K^{alg}/K_0)$ permute transitivement les racines primitives 
$p^r$-i\`emes de l'unit\'e, on en d\'eduit le r\'esultat annonc\'e.
\end{proof}
\par
Par ailleurs, nous aurons besoin du lemme ci-dessous, qui \'etablit la 
lin\'earisabilit\'e 
locale de l'action d'un automorphisme d'ordre une puissance de $p$ 
au voisinage d'un point fixe.
\begin{lem}
Soit $\sigma$ un $R$-automorphisme d'ordre $p^r$ de $R[[Z]]$, 
fixant $0$. Il existe $\rho$ 
dans $K$, avec $v_K(\rho)>0$, une racine primitive $p^r$-i\`eme de 
l'unit\'e $\zeta_{(r)}^l$ et un param\`etre $Z'$ de 
$R\{\displaystyle{\frac{Z}{\rho}}\}$ 
tels que $\sigma(Z')=\zeta_{(r)}^lZ'$. Autrement dit, $\sigma$ est 
lin\'earisable sur un sous-disque ferm\'e de centre $0$.
\end{lem}
\begin{proof}
Notons $B:=R[[Z]]$, $A:=R[[Z]]^{\sigma}=R[[T]]$. L'extension de corps 
$Fr B/Fr A$ est donn\'ee, d'apr\`es la th\'eorie de Kummer, par une 
\'equation $Y^{p^r}=u$, o\`u $u \in A$, $\bar u \ne 0$. Si on restreint 
le disque 
\`a un sous-disque ferm\'e (de param\`etre 
$Z_0=\displaystyle{\frac {Z}{\rho}}$) centr\'e 
en $0$, on aura $R\{Z_0\}^{\sigma}=R\{T_0\}$ pour un $T_0$ convenable 
(par exemple la norme de $Z_0$). Si la valuation de $\rho$ est suffisamment 
grande, $0$ est le seul point de ramification, on peut alors supposer que 
$u$ s'\'ecrit $u=T_0^h(1+\sum_{\nu>0}a_\nu T_0^{\nu})$, 
avec $v_K(a_{\nu})>0$ pour tout $\nu>0$ et $(h,p)=1$. En utilisant Bezout, 
on peut supposer $h=1$, et finalement, quitte \`a changer le param\`etre 
$T_0$, que $u=T_0$. Mais $\frac{R\{T_0\}[Y]}{(Y^{p^r}-T_0)}$ est 
int\'egralement clos et contenu dans $R\{Z_0\}$, donc finalement, 
$R\{Z_0\}=\frac{R\{T_0\}[Y]}{(Y^{p^r}-T_0)}$ et $Z':=Y$ est un 
param\`etre de $R\{Z_0\}$ qui convient.
\end{proof}
\begin{rem}
\label{linearisable}
La preuve ci-dessus montre que si $r=1$ et si $\sigma$ 
poss\`ede un unique point fixe, l'automorphisme 
$\sigma$ est lin\'earisable sur $R[[Z]]$, on retrouve ainsi la proposition 
6.2.1 de \cite{G-M 2}. C'est une question ouverte de savoir si, pour $r>1$, le 
r\'esultat est encore vrai : La remarque 6.2.2 de \cite{G-M 2} semble 
erron\'ee. Toutefois, si 
$\sigma^{p^{r-1}}$ poss\`ede un unique point fixe, alors $\sigma$ est encore 
lin\'earisable sur $R[[Z]]$.
\end{rem}
\subsection{Rel\`evement local cyclique pour les points doubles}
\begin{prop}
\label{relloc}
(i) Pour tout entier $r \ge 1$, quitte \`a faire une 
extension finie de $K$, il existe un entier $e \ge 1$ et un automorphisme 
d'ordre $p^r$ de $\mathcal C_e$ qui induit un automorphisme d'ordre $p^r$ 
de chaque branche de $\mathcal C_{e,k}$.\\
(ii) Si $\bar{\sigma}$ est un automorphisme de $\spec 
\displaystyle{\frac{k[[z_1,z_2]]}{(z_1z_2)}}$, ne permutant pas les branches 
du point double, qui induit un automorphisme d'ordre $p^r$ ($1 \le r \le 2$) 
de chaque branche, alors, quitte \`a faire une extension finie de $K$, il 
existe un entier $e \ge 1$ et un $R$-automorphisme $\sigma$ d'ordre 
$p^r$ de $\mathcal C_e$ qui induit $\bar{\sigma}$ sur 
$\mathcal C_{e,k}= \spec \displaystyle{\frac{k[[z_1,z_2]]}{(z_1z_2)}}$.
\end{prop}
\begin{figure}[htpb]
         \input{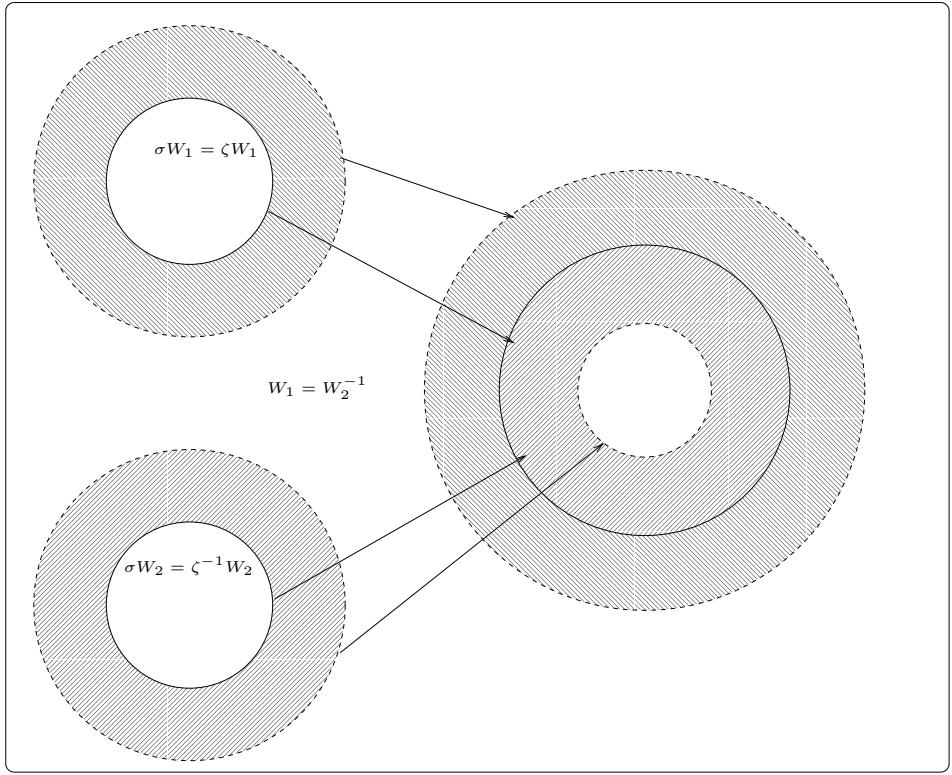}
         \caption{Construction d'un automorphisme de couronnes par recollement}
\end{figure}
\begin{proof}
Soit $r \ge 1$, on fixe une racine primitive $p^r$-i\`eme de l'unit\'e 
$\zeta_{(r)}$ dans $K^{alg}$.
On consid\`ere deux $R$-automorphismes $\sigma_1$ et $\sigma_2$ d'ordre $p^r$ 
de $\mathcal D$ fixant $0$, on supposera les points de $F_{\sigma_1}$ 
(resp. $F_{\sigma_2}$) rationnels sur $K$. Pour $i=1,2$, il existe un 
unique $h_i$ inversible dans $(\mathbb{Z}/p^r\mathbb{Z})$ tel que 
$\displaystyle{\frac{d\sigma_i}{dZ}}(0)=\zeta_{(r)}^{h_i^{-1}}$. 
D'apr\`es le lemme qui pr\'ec\`ede, il existe  
un param\`etre $W_i$ d'un sous-disque ferm\'e 
$D_i':=\{\omega \in D|v_K(\omega) \ge e_i\}$ de $D:=\mathcal D_K$ 
centr\'e en $0$ tel que $\sigma_i(W_i)=\zeta_{(r)}^{h_i^{-1}}W_i$. 
En utilisant le lemme \ref{ajustement}, on peut supposer que 
$h_1=-h_2$ (noter que $\bar \sigma_2$ ne change pas). On construit alors un 
automorphisme d'une couronne formelle de la fa\c con suivante :
\par
Soit $e_0$ un entier positif ou nul, notons $\sigma_0$ l'automorphisme de 
$\displaystyle{\frac{R\{W_1,W_2\}}{(W_1W_2-\pi^{e_0})}}$ d\'efini par 
$\sigma_0(W_1)=\zeta_{(r)}^{h_1^{-1}}W_1$, 
$\sigma_0(W_2)=\zeta_{(r)}^{h_2^{-1}}W_2$. Notons, pour $i=1,2$, $j_i$ 
l'inclusion 
$\displaystyle{\frac{R[[Z_i]]\{X_i\}}{(Z_iX_i-\pi^{e_i})}} \rightarrow 
R\{X_i,X_i^{-1}\}=R\{W_i,W_i^{-1}\}$ (o\`u $X_i$ s'\'ecrit 
$X_i=W_i^{-1}(1+\pi g_i)$ dans $R\{X_i,X_i^{-1}\}$). On munit, 
pour $i=1,2$, $\displaystyle{\frac{R[[Z_i]]\{X_i\}}{(Z_iX_i-\pi^{e_i})}}$ 
(resp. $\displaystyle{\frac{R\{W_1,W_2\}}{(W_1W_2-\pi^{e_0})}}$), 
d'une structure de $R[\mathbb Z/p^r\mathbb Z]$-module par 
$l.f_i:=\sigma_i^lf_i$ (resp. $l.f_0:=\sigma_i^lf_0$) pour 
$l\in \mathbb Z/p^r\mathbb Z$. 
On d\'efinit alors un morphisme 
$$\theta:\frac{R[[Z_1]]\{X_1\}}{(Z_1X_1-\pi^{e_1})} \times 
\frac{R[[Z_2]]\{X_2\}}{(Z_2X_2-\pi^{e_2})} \times 
\frac{R\{W_1,W_2\}}{(W_1W_2-\pi^{e_0})} \rightarrow 
R\{W_1,W_1^{-1}\} \times
R\{W_2,W_2^{-1}\}$$ 
de $R[\mathbb Z/p^r\mathbb Z]$-modules par 
$\theta(f_1,f_2,f_0)=(j_1(f_1)-f_0,j_2(f_2)-f_0)$. On obtient ainsi une 
structure de $R[\mathbb Z/p^r\mathbb Z]$-module sur le noyau de $\theta$, 
qui s'identifie au $R$-module 
$\displaystyle{\frac{R[[Z'_1,Z'_2]]}{(Z'_1Z'_2-\pi^{e_0+e_1+e_2})}}$. 
(En effet, on v\'erifie ais\'ement que $\ker \bar \theta$ est \'egal \`a 
$\displaystyle{\frac{k[[z_1,z_2]]}{(z_1z_2)}}$, et donc par le lemme \ref{lat} 
(ii), on montre que $\ker \theta$ est isomorphe \`a 
$\kouro$, avec $e=e_0+e_1+e_2$). On en d\'eduit un 
automorphisme d'ordre $p^r$ d'une couronne formelle d'\'epaisseur 
$e$ qui induit $\bar \sigma$.
\par
L'assertion (ii) r\'esulte alors de \cite{G-M 1}, th\'eor\`emes II 4.1 et 5.5. 
L'assertion (i) r\'esulte quant \`a elle de \cite{G-M 2} II, paragraphe 3.3.3. 
\end{proof}

\begin{cor}
\label{relevlocalptdb}
On consid\`ere des entiers $a$ et $n$, avec $1\le a \le 2$ et $(n,p)=1$. 
Soit $\bar \sigma$ un automorphisme d'ordre $p^an$ de 
$\displaystyle{\frac{k[[z_1,z_2]]}{(z_1z_2)}}$ , ne 
permutant pas les branches, induisant un automorphisme d'ordre $p^an$ de 
chaque branche. On suppose de plus que $\bar \sigma^{p^a}$ agit de mani\`ere 
kumm\'erienne, c'est \`a dire que quitte \`a changer de param\`etres 
$(z_1,z_2)$, $\bar \sigma^{p^a}(z_1)=\bar \theta z_1$ et 
$\bar \sigma^{p^a}(z_2)=\bar \theta^{-1} z_2$, o\`u 
$\bar \theta$ est une racine 
primitive $n$-i\`eme de l'unit\'e dans $k$. 
Alors, quitte \`a faire une extension 
finie de $K$, $\bar \sigma$ se rel\`eve en un automorphisme d'ordre 
$p^an$ de $\mathcal C_e$, pour un entier $e \ge 1$ convenable.
\end{cor}
\begin{proof}
\label{relcycptdb}
Comme $\bar \sigma^{p^a}$ agit de mani\`ere kumm\'erienne, quitte \`a changer 
de param\`etres, on peut supposer $\bar \sigma^{p^a}(z_1)=\bar \theta z_1$ et 
$\bar \sigma^{p^a}(z_2)=\bar \theta^{-1} z_2$, o\`u $\theta$ est une racine 
primitive $n$-i\`eme de l'unit\'e dans $R$. Le sous-anneau de 
$\displaystyle{\frac{k[[z_1,z_2]]}{(z_1z_2)}}$ form\'e des \'el\'ements fixes 
sous $\bar \sigma^{p^a}$ est alors 
$\displaystyle{\frac{k[[x_1,x_2]]}{(x_1x_2)}}$, o\`u $x_i=z_i^n$. 
L'automorphisme $\bar \sigma^n$ induit un automorphisme $\bar \tau$ d'ordre 
$p^a$ de $\displaystyle{\frac{k[[x_1,x_2]]}{(x_1x_2)}}$. Cet automorphisme se 
rel\`eve alors en un automorphisme $\tau$ d'ordre $p^a$ de 
$\displaystyle{\frac{R[[X_1,X_2]]}{(X_1X_2-\pi^{en})}}$ , pour 
un entier $e$ convenable (quitte \`a faire une extension finie de $K$). 
Ecrivons 
$$\tau(X_1)=X_1(1+\pi b+\sum_{\nu>0}(a_{\nu}X_1^{\nu}+a_{-\nu}X_2^{-\nu})).$$ 
Consid\'erons la $R$-alg\`ebre 
$\displaystyle{\frac{R[[Z_1,Z_2]]}{(Z_1Z_2-\pi^{e})}}$, o\`u $Z_i$ rel\`eve 
$z_i$ pour $i=1,2$. En identifiant $X_i=Z_i^n$, 
$\displaystyle{\frac{R[[X_1,X_2]]}{(X_1X_2-\pi^{en})}}$ est un sous-anneau de 
$\displaystyle{\frac{R[[Z_1,Z_2]]}{(Z_1Z_2-\pi^{e})}}$. 
Posons 
$$\tilde \tau(Z_1)=Z_1(1+\pi b+\sum_{\nu>0}(a_{\nu}Z_1^{n\nu}+a_{-\nu}Z_2^{-n\nu}))^{\frac 1n},$$ $\tilde \tau$ prolonge alors $\tau$ et est d'ordre 
$p^a$. De plus $\tilde \tau$ commute avec l'automorphisme $\mu$ d'ordre $n$ 
d\'efini par $\mu(Z_1)=\theta Z_1$. Par suite, le groupe engendr\'e par 
$\tilde \tau$ et $\mu$ est cyclique d'ordre $p^an$, et un g\'en\'erateur convenable 
rel\`eve $\bar \sigma$.  
\end{proof}
\subsection{Rel\`evement local pour le groupe di\'edral}
Soit $D$ le groupe di\'edral d'ordre $2np^r$, o\`u $n$ est un entier 
premier \`a $p$ et $r$ un entier inf\'erieur ou \'egal \`a $2$. Le groupe $D$ admet une pr\'esentation 
$$D=\langle \sigma,\tau|\sigma^{np^r}=1,\tau^2=1,
\tau \sigma \tau= \sigma^{-1}\rangle.$$
Les m\'ethodes formelles ci-dessus permettent encore de construire un 
rel\`evement d'une action de $D$ sur un point double :
\begin{prop}
\label{reldihptdb}
On consid\`ere une action de $D$ sur 
$\spec \displaystyle{\frac{k[[z_1,z_2]]}{(z_1z_2)}}$, 
telle que $\tau$ permute les branches du point double, $\sigma$ les fixe, 
et $\sigma$ induit un automorphisme d'ordre $np^r$ de chaque branche. Alors, 
l'action de $D$ se rel\`eve en une action de $D$ sur une couronne formelle 
sur $R:=W(k)[\zeta_{(r)}]$ d'\'epaisseur paire, où $\zeta_{(r)}$ est une 
racine primitive $p^r$-ième de l'unité.
\end{prop}
\begin{proof}
L'action de $D$ peut s'\'ecrire, pour un bon choix de coordonn\'ees 
$(z_1,z_2)$, $\tau z_1=z_2$, $\sigma(z_1)=f(z_1)$, o\`u 
$f(x) \in k[[x]]$ est une s\'erie qui d\'efinit un automorphisme 
d'ordre $np^r$ de $k[[x]]$. Remarquons que cette \'ecriture et la relation 
$\tau \sigma \tau= \sigma^{-1}$ entra\^ \i nent que l'action du groupe 
cyclique engendr\'e par $\sigma$ est kumm\'erienne. 
On peut relever la s\'erie $f$ en une s\'erie qui d\'efinit un automorphisme 
d'ordre $np^r$ de $R[[X]]$, encore not\'e $\sigma$ et qui fixe $0$. 
Il existe alors un sous-disque ferm\'e formel 
$\spec R\{X_0\}$ de $\spec R[[X]]$, centr\'e en $0$, fix\'e par $\sigma$, 
tel que $\sigma(X_0)=\mu X_0$, o\`u $\mu$ est une racine primitive 
$np^r$-i\`eme de l'unit\'e dans $R$. Notons $e$ l'\'epaisseur de la couronne 
formelle ``compl\'ementaire''. 
On en d\'eduit de mani\`ere analogue \`a la preuve de \ref{relcycptdb}
une suite exacte de $R$-modules
$$0\to \frac{R[[Z_1,Z_2]]}{(Z_1Z_2-\pi^{2e})}\to 
R[[X]]\{\frac{\pi^e}{X}\} \times R[[Y]]\{\frac{\pi^e}{Y}\}
\to R\{X_0,X_0^{-1}\} \to 0
$$
En fait, le morphisme surjectif est \'equivariant sous l'action de $D$, 
si on munit le $R$-module du milieu de l'action $\sigma(g(X),h(Y)):=
(g(\sigma(X)),h(\sigma^{-1}(Y))$ et $\tau(X)=Y$. On obtient ainsi une action 
de $D$ sur le disque formel (d'\'epaisseur $2e$) 
$$\spec \displaystyle{\frac{R[[Z_1,Z_2]]}{(Z_1Z_2-\pi^{2e})}}$$
qui rel\`eve l'action de $D$ sur le point double.
\end{proof}
\subsection{Application au rel\`evement galoisien}
On consid\`ere une courbe alg\'ebrique $Y$ sur $k$, propre, connexe, 
nodale, muni d'un sous-groupe fini $G$ de $\aut_k Y$, op\'erant librement
 sur un ouvert dense. Soit $y$ un point double de $Y$, de stabilisateur $I_y$ 
cyclique. Le sous-groupe $H_y$ de $I_y$ form\'e des automorphismes ne 
permutant pas les branches analytiques du point double $y$ et d'ordre premier 
\`a $p$ est muni de deux caract\`eres $\chi_{y,1},\chi_{y,2}:H_y
\rightarrow k^{*}$ correspondant \`a l'action de $H_y$ sur l'espace 
tangent en $y$. Suivant \cite{Sa1}, on dira que l'action de $G$ sur $Y$ est 
kumm\'erienne si, pour tout point double $y$ de $Y$, de stabilisateur cyclique,
 on a $\chi_{y,1}\chi_{y,2}=1$. Cette d\'efinition co\" \i ncide avec la 
pr\'ec\'edente 
lorsque $G$ est cyclique d'ordre premier \`a $p$.
\begin{thm}
\label{relgal}
Supposons l'action de $G$ kumm\'erienne et que pour un point 
ferm\'e $y$ de $Y_0$ :\\
1. Si y est un point lisse, le groupe d'inertie $I_y$ de $y$ est 
cyclique d'ordre 
$n(y)p^{r(y)}$, avec $(n(y),p)=1$ et $0\le r(y)\le 2$.\\
2. Si y est un point double, on a alors (i) ou bien (ii) :
\begin{itemize}
\item
(i) Le groupe d'inertie $I_y$ de $y$ est cyclique d'ordre $n(y)p^{r(y)}$ 
et l'action de $I_y$ sur les branches est triviale.
\item
(ii) Le groupe d'inertie $I_y$ de $y$ est di\'edral d'ordre $2n(y)p^{r(y)}$, 
avec $(n(y),p)=1$ et $0\le r(y)\le 2$, de pr\'esentation 
$$I_y=\langle \sigma,\tau|\sigma^{n(y)p^{r(y)}}=1,\tau^2=1,
\tau \sigma \tau= \sigma^{-1}\rangle.$$
Les branches du point double $y$ sont permut\'ees par $\tau$ et $\sigma$ 
induit un automorphisme d'ordre $n(y)p^{r(y)}$ de chacune des deux branches.
\end{itemize}
 Quitte \`a faire une extension finie de $K$, il existe un 
$R$-sch\'ema $\mathcal Y$ normal, propre, plat, de 
fibre sp\'eciale $Y$, de fibre g\'en\'erique lisse et g\'eom\'etriquement 
connexe sur 
$K$, muni d'une action de $G$ relevant celle sur $Y$.
\end{thm}
\begin{proof}\ \\
\par
Soit $X:=Y/G$ le quotient de $Y$ par $G$. C'est encore une courbe nodale
 (voir \cite{dJ} prop. 4.2). Notons $B$ l'ensemble des points de branchement 
dans $X$ du morphisme quotient $\overline{f}:Y\longrightarrow X$, et 
$Ram$ l'ensemble des points de ramification de 
$\overline{f}$. On choisit, pour $x \in B$, un point $y_x$ de $Ram$ tel que 
$f(y_x)=x$. Soit $B_1$ (resp. $B_2$, $B_3$ ) le sous-ensemble de $X$ 
form\'e des $x$ tels que $y_x$ soit un point lisse de $Y$ (resp. un point 
double d'inertie cyclique, un point double d'inertie di\'edrale). Remarquons 
que pour $x$ dans $B_1 \cup B_3$ (resp. $B_2$), $x$ est un point lisse de $X$ 
(resp. $x$ est un point double de $X$).
\par
Pour $x$ dans $B$, il existe un rel\`evement local 
$$
  \diagram
     I_{y_x} \rto \drto & \aut_R \mathcal A_x \dto \\
                    & \aut_k \mathcal O_{Y,y_x}
  \enddiagram
$$
o\`u $\mathcal A_x$ est isomorphe \`a $R[[Z]]$ si $x$ appartient \`a $B_1$ et 
$\mathcal A_x$ est isomorphe \`a 
$\displaystyle{\frac{R[[Z_1,Z_2]]}{(Z_1Z_2-\pi^{e_x})}}$ pour $x$ dans 
$B_2 \cup B_3$, o\`u $e_x$ est un entier sup\'erieur ou \'egal \`a $1$, qui 
est de plus pair pour $x$ dans $B_2$.
\par
Pour $x$ un point de $B$, on note $U_x$ un voisinage affine de $x$ dans $X$, 
inclus dans la r\'eunion des composantes irr\'eductibles de $X$ qui 
contiennent $x$, et $V_x:=\overline{f}^{-1}(U_x)$. On note $U$ l'ouvert dense 
$X\setminus B$ et $V:=\overline{f}^{-1}(U)$.
\par
D'apr\`es \ref{modele}, il existe un $R$-sch\'ema $\mathcal X$ 
propre et plat, normal, de fibre 
g\'en\'erique lisse sur $K$, de fibre sp\'eciale $X$, dont l'\'epaisseur au 
point double $x$ de $B_2$ est $pe_x$. On notera $\hat{\mathcal X}$ 
le $R$-sch\'ema formel compl\'et\'e 
de $\mathcal X$ le long de sa fibre sp\'eciale, et $\mathcal{U}_x$ l'ouvert 
de $\hat{\mathcal X}$ correspondant \`a $U_x$, pour $x$ dans $B$. 
\par
Pour $x$ dans $B$, il existe un morphisme \'etale $U_x' \rightarrow U_x$ 
(qu'on peut supposer surjectif, quitte \`a restreindre $U_x$) tel que le 
morphisme 
$V_x' \rightarrow U_x'$ obtenu par changement de base \`a partir de 
$V_x \rightarrow U_x$ soit d\'ecompos\'e, ie 
$V_x'=\ind_{I_{y_x}}^G W_x'$, o\`u 
$W_x'$ est la composante connexe de $V_x'$ contenant l'unique point au-dessus
 de $y_x$. En utilisant les rel\`evements locaux choisis et 
le th\'eor\`eme \ref{recolloc}, on construit un  rel\`evement 
$\mathcal{W}_x' \rightarrow \mathcal{U}_x'$ de $\overline{f}_{|W_x'}$ 
compatible avec l'action de $I_{y_x}$.
\par
Posons alors $\mathcal{V}_x':=\ind_{I_{y_x}}^G \mathcal{W}_x'$. C'est un 
rel\`evement de $\overline{f}_{|V_x'}$ compatible avec l'action de $G$. 
Un argument de descente \'etale (voir \cite{Ga} lemme 3.6) permet alors 
de construire un rel\`evement $\mathcal{V}_x \rightarrow \mathcal{U}_x$, 
galoisien de groupe $G$, de $V_x \rightarrow U_x$.  L'unicit\'e du 
rel\`evement du lieu \'etale nous permet de recoller ces rel\`evements en un 
rel\`evement $\hat{\mathcal Y} \stackrel{\hat{f}}{\rightarrow} \hat{\mathcal X}$ de $\overline{f}$, galoisien de groupe $G$. 
Comme $\hat f$ est propre, le th\'eor\`eme d'alg\'ebrisation des faisceaux 
coh\'erents (\cite{EGA} corollaire 5.1.6) permet de 
construire une $\mathcal O_{\mathcal X}$-alg\`ebre  cohérente $\mathcal A$ 
dont le complété s'identifie \`a la  $\mathcal O_{\hat{\mathcal X}}$-alg\`ebre cohérente $\mathcal O_{\hat{\mathcal Y}}$. Posons alors 
$\mathcal Y:=\spec \mathcal A$. L'action du groupe $G$ sur
$\hat{\mathcal Y}$ se relève de manière unique en une action de $G$ sur 
$\mathcal Y$. Le couple $(\mathcal Y,G)$ convient.
\end{proof}


\begin{thebibliography}{17}
\bibitem{B-M} Bertin J., M\'ezard A. : D\'eformations formelles des 
rev\^etements sauvagement ramifi\'es de courbes alg\'ebriques. 
Pr\'epublication $n^o 439$ de l'institut Fourier, (1998)
\bibitem{Ga} Garuti M. : Prolongement de rev\^etements galoisiens en 
g\'eom\'etrie 
rigide. Compositio Math. 104 ($n^o 3$), 305-331 (1996)
\bibitem{G-M 1} Green B., Matignon M. : Lifting of Galois covers of smooth 
curves. Compositio Math. 113 ($n^o$ 3), 237-272 (1998)
\bibitem{G-M 2} Green B., Matignon M. : Order $p$ automorphisms of the open 
disc of a $p$-adic field. J. Amer. Soc 12 ($n^o$ 1), 269-303 (1999)
\bibitem{EGA} Grothendieck A., Dieudonné J. : Eléments de G\'eom\'etrie 
Alg\'ebrique, chap. III, Publ. Math. IHES 11 (1961)
\bibitem{SGA1} Grothendieck A. : SGA1, Rev\^etements \'etales et 
groupe fondamental (S\'eminaire de g\'eom\'etrie alg\'ebrique du Bois-Marie), 
LNM 224 Springer-Verlag (1960-61)
\bibitem{H} Harbater D. : Formal patching and adding branch points. Amer. J. 
of Math. 115 ($n^o 3$), 487-508 (1993)
\bibitem{H-S} Harbater D., Stevenson K. : Patching and thickening problems. 
J. of Alg. 212, 272-304 (1999)
\bibitem{Heth} Henrio Y. : Arbres de Hurwitz et automorphismes d'ordre $p$ 
des disques et des couronnes $p$-adiques formels. Th\`ese, universit\'e 
Bordeaux I (1999)
\bibitem{He} Henrio Y. : Automorphismes d'ordre $p$ des couronnes $p$-adiques 
ouvertes. C. R. Acad. Sci Paris, t. 329, S\'erie I, 47-50 (1999)
\bibitem{dJ} de Jong A.J. : Families of curves and alterations. Ann. Inst. 
Fourier, Grenoble 47 ($n^o 2$), 599-621 (1997)
\bibitem{O1} Oort F. : Lifting algebraic curves, abelian varieties, 
and their endomorphisms to characteristic zero. Proceedings of Symposia 
in Pure Mathematics, Vol. 46 (1987)
\bibitem{O2} Oort F. : Some Questions in Algebraic Geometry. Utrecht Univ., 
Math. Dept. Preprint Series, June 1995.
\bibitem{OSS} Oort F., Sekiguchi T., Suwa N. : On the deformation of 
Artin-Schreier to Kummer. Ann. scient. Ec. Norm. Sup. 22, 345-375 (1989)
\bibitem{Ra} Raynaud M. : Rev\^etements de la droite affine en 
caract\'eristique $p>0$ et conjecture d'Abhyankar. Invent. Math 116, 
425-462 (1994)
\bibitem{Sa1} Sa\" idi M. : Rev\^etements \'etales ab\'eliens, courants sur 
les graphes et r\'eduction semi-stable des courbes. Manuscripta Math. 89, 
245-265 (1996)
\end{thebibliography}
\end{document}